\documentclass[preprint]{elsarticle}

\usepackage{amssymb}
\usepackage{amsfonts}
\usepackage{amsmath}
\usepackage{graphicx}
\usepackage{amsthm}
\usepackage{algorithm}
\usepackage{epic}

\newtheorem{rem}{Remark}

\newtheorem{exmpl}{Example}

\newtheorem{lemma}{Lemma}

\begin{document}

\begin{frontmatter}

\title{Fast multiplication, determinants, and inverses of arrowhead and diagonal-plus-rank-one matrices over associative fields} 
\author[FESB]{N.\ Jakov\v{c}evi\'{c} Stor\corref{c1}}
\ead{nevena@fesb.hr}
\author[FESB]{I.\ Slapni\v{c}ar}
\ead{ivan.slapnicar@fesb.hr}

\cortext[c1]{Corresponding author}

\address[FESB]{Faculty of Electrical Engineering, Mechanical Engineering and
Naval Architecture, University of Split, Rudjera Bo\v{s}kovi\'{c}a
32, 21000 Split, Croatia}

\begin{abstract}
The article considers arrowhead and diagonal-plus-rank-one matrices  in  $\mathbb{F}^{n\times n}$ where $\mathbb{F}\in\{\mathbb{R},\mathbb{C},\mathbb{H}\}$.
$\mathbb{H}$  is a non-commutative field of quaternions. We give unified formulas for fast matrix-vector multiplications, determinants, and inverses for considered matrices. The formulas are unified in the sense that the same formula holds in both, commutative and noncommutative algebras. Each formula requires $O(n)$ arithmetic operations. Most of the formulas hold for block matrices, as well.
\end{abstract}

\begin{highlights}
\item Fast $O(n)$ formulas.
\item Matrix elements are real or complex numbers, quaternions, or block matrices.
\item Same formulas for all types of matrix entries.
\item Code and examples on GitHub.
\end{highlights}


\end{frontmatter}

\section{Introduction and definitions}

Arrowhead matrices and diagonal-plus-rank-one matrices arise in many applications and computations with such matrices are parts of many important linear algebra algorithms (for details see \cite{GV96, NEG14, JSB15}). We give unified formulas for matrix-vector multiplications, determinants, and inverses for both types of matrices having elements in commutative and noncommutative associative algebras. Our results complement and extend the existing results in the literature. 

All matrices are in $\mathbb{F}^{n\times n}$ where $\mathbb{F}\in\{\mathbb{R},\mathbb{C},\mathbb{H}\}$, where $\mathbb{H}$ is a non-commutative field of quaternions. Each formula requires $O(n)$ arithmetic operations.

In Section 1, we give basic definitions. In section 2, we give formulas for the multiplication of vectors with both types of matrices. In Section 3, we give formulas for determinants, and in Section 4, we give formulas for inverses of both types of matrices. In Section 5, we explain how the formulas apply to block matrices. Discussion and conclusions are in Section 6.

\subsection{Quaternions}

{\em Quaternions} are a non-commutative associative number system that extends complex numbers, introduced by Hamilton \cite{Ham53, Ham66}. For {\em basic quaternions} $\mathbf {i}$, 
$\mathbf {j}$, and $\mathbf {k}$, the quaternions have the form
$$q=a+b\ \mathbf {i} +c\ \mathbf {j} +d\ \mathbf {k},\quad a,b,c,d, \in \mathbb{R}.$$
The {\em multiplication table} of basic quaternions is the following:
$$\begin{array}{c|ccc}
\times & \mathbf {i} & \mathbf {j} & \mathbf {k} \\ \hline
\mathbf {i} & -1 & \mathbf {k} & -\mathbf {j} \\
\mathbf {j} & -\mathbf {k} & -1 & \mathbf {i} \\
\mathbf {k} & \mathbf {j} & -\mathbf {i} & -1
\end{array}$$
{\em Conjugation} is given by 
$$\bar q=a-b\ \mathbf {i} -c\ \mathbf {j} -d\ \mathbf {k}.$$
Then, 
$$\bar q q=q\bar q=|q|^2=a^2+b^2+c^2+d^2.$$
%
%
%
Let $f(x)$ be a complex analytic function. The value $f(q)$, where $q\in\mathbb{H}$, is computed by evaluating the extension of $f$ to the quaternions at $q$ (see \cite{Sud79}). 
%
%
All of the above is implemented in the Julia \cite{Julia} package Quaternions.jl \cite{Quaternions}.

%

Quaternions are {\em homomorphic} to $\mathbb{C}^{2\times 2}$:
$$
q\to \begin{bmatrix}a+b\,\mathbf{i} & c+d\, \mathbf{i}\\-c+d\, \mathbf{i} & a-b\, \mathbf{i}\end{bmatrix}\equiv C(q),$$
with eigenvalues $q_s$ and $\bar q_s$.

\subsection{Arrowhead and DPR1 matrices}

Let $\star$ denote the transpose of a real matrix, the conjugate transpose (adjoint) of a complex or quaternionic matrix, and the conjugate of a scalar.

The {\em arrowhead matrix} (Arrow) is a matrix of the form
$$
A=\begin{bmatrix} D & u \\v^\star & \alpha \end{bmatrix}\equiv \operatorname{Arrow}(D,u,v,\alpha),$$
where  $D$ is a diagonal matrix with diagonal elements $d_i\equiv D_{ii}$, $\operatorname{diag}(D),u, v \in\mathbb{F}^{n-1}$, and $\alpha\in\mathbb{F}$, or any symmetric permutation of such matrix. 

The {\em diagonal-plus-rank-one matrix} (DPR1) is a matrix of the form 
$$A=\Delta+x \rho y^\star\equiv \operatorname{DPR1}(\Delta,x,y,\rho),$$
where  $\Delta$ is diagonal matrix with diagonal elements $\delta_i=\Delta_{ii}$, $\operatorname{diag}(\Delta),x, y \in\mathbb{F}^{n}$, and $\rho \in \mathbb{F}$. 

\section{Matrix-vector multiplication}

\begin{lemma} \label{lemma0} 
Let $A=\operatorname{Arrow}(D,u,v,\alpha)$ be an arrowhead matrix with the tip at position $A_{ii}=\alpha$, and let $z$ be a vector. Then $w=Az$, where
$$
\begin{aligned}
w_j&=d_jz_j+u_jz_i, \quad i=1,2,\cdots,i-1\\
w_i&=v_{1:i-1}^\star z_{1:i-1} +\alpha z_i + v_{i:n-1}^\star z_{i+1:n} \\
w_j&=u_{j-1}z_i+d_{j-1}z_j,\quad j=i+1,i+2,\cdots,n.
\end{aligned}$$
Further, let $A=\operatorname{DPR1}(\Delta,x,y,\rho)$ be a DPR1 matrix and let $\beta=\rho(y^\star x) \equiv \rho (y\cdot x)$. Then $w=Az$, where
$$
w_i=\delta_i z_i+x_i\beta,\quad i=1,2,\cdots,n.$$
\end{lemma}

\proof The formulas follow by direct multiplication. \qed

%

\section{Determinants}

Determinants are computed using two basic facts:
the determinant of the triangular matrix is a product of diagonal elements ordered from the first to the last, and
the determinant of the product is the product of determinants.

If $\mathbb{F}\in \{\mathbb{R},\mathbb{C}\}$, we have the following lemmas.

%
%

\begin{lemma} \label{lemma1} 
Let $A=\operatorname{Arrow}(D,u,v,\alpha)$ be an arrowhead matrix. If all $d_i\neq 0$, the determinant of $A$ is equal to
\begin{equation}
\det(A)= (\prod_i d_i)(\alpha-v^\star D^{-1}u). \label{da1}
\end{equation}
If $d_i=0$, then 
\begin{equation} \label{da2}
\det(A)=-\bigg(\prod_{j=1}^{i-1}d_j\bigg)\cdot  v_i^\star 
\cdot \bigg(\prod_{j=i+1}^{n-1}d_j\bigg)\cdot  u_i.
\end{equation}
\end{lemma}

\proof The proof is modeled after Proposition 2.8.3, Fact 2.14.2 and Fact 2.16.2 from \cite{Ber09}.
The formula (\ref{da1}) follows from the factorization
$$
\begin{bmatrix} D & u  \\ v^\star & \alpha \end{bmatrix} 
=\begin{bmatrix} D & 0 \\ 0  & 1 \end{bmatrix}
\begin{bmatrix} I & 0 \\ v^\star & 1 \end{bmatrix}
\begin{bmatrix} I & D^{-1} u \\ 0 & \alpha - v^\star D^{-1} u \end{bmatrix}.
$$
The formula (\ref{da2}) is the only non-zero term in the standard expansion of the determinant of the matrix
$$
D=\begin{bmatrix}D_1 & 0 & 0 & u_1 \\ 0 & 0 & 0 & u_i \\ 0 & 0 & D_2 & u_2 \\
v_1^\star & v_i^\star & v_2^\star & \alpha \end{bmatrix}.$$ 
\qed 

\begin{lemma} \label{lemma2} 
Let $A=\operatorname{DPR1}(\Delta,x,y,\rho)$ be a DPR1 matrix. If all $\delta_i\neq 0$, the determinant of $A$ is equal to
\begin{equation} \label{dd1}
\det(A)=(\prod_i \delta_i)(1+y^\star \Delta^{-1}x\rho).
\end{equation}
If $\delta_i=0$, then 
\begin{equation} \label{dd2}
\det(A)=(\prod_{j=1}^{i-1} \delta_j) y_i^\star (\prod_{j=i+1}^n \delta_j)x_i \rho.
\end{equation}
\end{lemma}

\proof The proof is modeled after Fact 12.16.3 and Fact 12.16.4 from \cite{Ber09}.
The formula (\ref{dd1}) follows from the factorizations
$$
\Delta+x\rho y^\star=\Delta (I+\Delta^{-1}x\rho y^\star)$$
and
$$
\begin{bmatrix} I & 0\\ d^\star & 1 \end{bmatrix} 
\begin{bmatrix} I+cd^\star & c \\ 0 & 1\end{bmatrix} 
\begin{bmatrix} I & 0 \\ -d^\star & 1\end{bmatrix}= 
\begin{bmatrix} I & c \\ 0 & 1+d^\star c\end{bmatrix}.$$
The formula (\ref{dd2}) follows from Fact 2.14.2 from \cite{Ber09}.
\qed

\subsection{Matrices of quaternions}
The determinant of the matrix of quaternions can be defined using a determinant of its corresponding homomorphic complex matrix. For matrices of quaternions, the determinant is not well defined due to non-commutativity. However, the Study determinant, $|\det(A)|$, is well defined (see \cite{Asl96}). Therefore, after computing the respective determinant using the formulas from Lemma \ref{lemma1} and Lemma \ref{lemma2}, it suffices to compute the absolute value. 

\section{Inverses}


\begin{lemma} \label{lemma3}

Let $A=\operatorname{Arrow}(D,u,v,\alpha)$ be a non-singular arrowhead matrix with the 
tip at the position $A_{ii}=\alpha$ and let $P$ be the permutation matrix of the permutation $p=(1,2,\cdots,i-1,n,i,i+1,\cdots,n-1)$.
If all $d_j\neq 0$, the inverse of $A$ is a DPR1 matrix
\begin{equation} \label{ia1}
A^{-1} =\operatorname{DPR1}(\Delta,x,y,\rho)=\Delta+x\rho y^\star,
\end{equation}
where
$$
\Delta=P\begin{bmatrix}D^{-1} & 0\\ 0 & 0\end{bmatrix}P^T,
\ x=P\begin{bmatrix}D^{-1}u \\-1\end{bmatrix},\
y=P\begin{bmatrix}D^{-\star}v \\-1\end{bmatrix},\
\rho=(\alpha-v^\star D^{-1} u)^{-1}.
$$
If $d_j=0$, the inverse of $A$ is an arrowhead matrix with the tip of the arrow at position $(j,j)$ and zero at position $A_{ii}$ (the tip and the zero on the shaft change places). In particular, let $\hat P$ be the permutation matrix of the permutation
$\hat p=(1,2,\cdots,j-1,n,j,j+1,\cdots,n-1)$. Partition $D$, $u$ and $v$ as
$$
D=\begin{bmatrix}D_1 & 0 & 0 \\ 0 & 0 & 0 \\ 0 & 0 & D_2\end{bmatrix},\quad
u=\begin{bmatrix} u_1 \\ u_j \\u_2\end{bmatrix},\quad
v=\begin{bmatrix} v_1 \\ v_j \\v_2\end{bmatrix}.
$$
Then 
\begin{equation} \label{ia2} 
A^{-1}=P\begin{bmatrix} \hat D & \hat u\\ \hat v^\star & \hat \alpha \end{bmatrix}P^T,
\end{equation}
where
\begin{align*}
\hat D&=\begin{bmatrix}D_1^{-1} & 0 & 0 \\ 0 & D_2^{-1} & 0 \\ 0 & 0 & 0\end{bmatrix},\quad
\hat u= \begin{bmatrix}-D_1^{-1}u_1 \\ -D_2^{-1}u_2\\ 1 \end{bmatrix} u_j^{-1},\quad
\hat v= \begin{bmatrix}-D_1^{-\star}v_1 \\ -D_2^{-\star}v_2\\ 1\end{bmatrix}v_j^{-1},\\
\hat \alpha&=v_j^{-\star}\left(-\alpha +v_1^\star D_1^{-1} u_1+v_2^\star D_2^{-1}u_2\right) u_j^{-1}.
\end{align*}
\end{lemma}
\proof  The formula (\ref{ia1}) follows by multiplication. The formula (\ref{ia2}) is similar to the one from Section 2 of \cite{JSB15}. \qed

\begin{lemma} \label{lemma4}
Let $A=\operatorname{DPR1}(\Delta,x,y,\rho)$ be a non-singular DPR1 matrix. 
If all $\delta_j\neq 0$, the inverse of $A$ is a DPR1 matrix
\begin{equation} \label{id1} 
A^{-1} =\operatorname{DPR1}(\hat \Delta,\hat x,\hat y,\hat \rho)=\hat\Delta+\hat x\hat \rho \hat y^\star,
\end{equation}
where
$$
\hat \Delta=\Delta^{-1},\quad
\quad \hat x=\Delta^{-1}x,\quad
\hat y=\Delta^{-\star}y,\quad
\hat \rho=-\rho(I+y^\star \Delta^{-1} x\rho)^{-1}.
$$
If $\delta_j=0$, the inverse of $A$ is an arrowhead matrix with the tip of the arrow at position $(j,j)$.
In particular, let $P$ be the permutation matrix of the permutation
$p=(1,2,\cdots,j-1,n,j,j+1,\cdots,n-1)$.
Partition $\Delta$, $x$ and $y$ as
$$
\Delta=\begin{bmatrix}\Delta_1 & 0 & 0 \\ 0 & 0 & 0 \\ 0 & 0 & \Delta_2\end{bmatrix},\quad
x=\begin{bmatrix} x_1 \\ x_j \\x_2\end{bmatrix},\quad
y=\begin{bmatrix} y_1 \\ y_j \\y_2\end{bmatrix}.
$$
Then, 
\begin{equation} \label{id2} 
A^{-1} =P\begin{bmatrix} D & u \\ v^T & \alpha \end{bmatrix}P^T,
\end{equation}
where
\begin{align*}
D&=\begin{bmatrix} \Delta_1^{-1} & 0\\ 0 &\Delta_2^{-1}\end{bmatrix},\quad
u= \begin{bmatrix}-\Delta_1^{-1}x_1 \\ -\Delta_2^{-1}x_2\end{bmatrix} x_j^{-1},\quad
v= \begin{bmatrix}-\Delta_1^{-\star}y_1 \\ -\Delta_2^{-\star}y_2\end{bmatrix}y_j^{-1},\\
\alpha&=y_j^{-\star}\left(\rho^{-1} +y_1^\star \Delta_1^{-1} x_1+y_2^\star \Delta_2^{-1}x_2\right) x_j^{-1}.
\end{align*}
\end{lemma}
\proof The formula (\ref{id1}) follows form fact 2.16.3 of \cite{Ber09}. The formula (\ref{id2}) is similar to the one from Section 2 of \cite{JSB15b}. \qed

\section{Block matrices}

If elements of the matrices are themselves matrices in $\mathbb{F}^{k\times k}$ where 
$\mathbb{F}\in\{\mathbb{R},\mathbb{C},\mathbb{H}\}$, the above formulas can be used as follows.

If, in Lemma \ref{lemma0}, the elements of the vector $z$ satisfy $z_i\in\mathbb{F}^{k\times k}$, $i=1,\ldots,n$, the formulas hold directly and yield a block vector $w$.

The formulas (\ref{da1})-(\ref{dd2}) from Lemmas \ref{lemma1} and \ref{lemma2} return a matrix in $\mathbb{F}^{k\times k}$, so one more step is required -- computing the determinant of this matrix.

The formulas for inverses (\ref{ia1})-(\ref{id2}) from Lemmas \ref{lemma3} and \ref{lemma4} return the corresponding block matrices provided all inverses within formulas are well defined. It is clearly possible that a block arrowhead or DPR1 matrix
is non-singular even if some of the blocks involved in formulas (\ref{ia1})-(\ref{id2}) are singular. However, in such cases, the respective inverses do not have the structure required by the formulas, 
so the formulas cannot be applied.

\section{Discussion and conclusion}

We have derived formulas for basic matrix functions of arrowhead and diagonal-plus-rank-one matrices whose elements are real numbers, complex numbers, or quaternions. Each formula requires $O(n)$ arithmetic operations, so they are optimal. 
Each formula is unified in the sense that the same formula is used for any type of matrix element, including block matrices.
When the formulas are applied to matrices of quaternions or block matrices, due to non-commutativity, the operations must be executed exactly in the order specified.

Our results complement and extend the existing results in the literature, like \cite{NEG14, Ber09, GV96, Asl96}.

The code, written in the programing language Julia \cite{Julia}, along with examples, is available on GitHub \cite{Inverses}. The code relies on the Julia's polymorphism (or multiple dispatch) feature. 

\section*{Funding} This work has been fully supported by Croatian Science Foundation under the project IP-2020-02-2240. 


\end{document}